\documentclass[preprint,final,3p,times,12.8pt,times]{elsarticle}
\usepackage{mathrsfs}
\usepackage{array,multirow,bm}
\usepackage{siunitx} 
\usepackage{enumerate} 
\usepackage{graphicx}
\usepackage{amsthm,amsmath,amssymb}  
\usepackage[normalem]{ulem}
\usepackage{algorithm,algorithmicx,algpseudocode} 
\usepackage[labelfont=it,textfont=it,font=small]{caption}
\usepackage{hyperref}

\newtheorem{theorem}{Theorem}[section]

\newtheorem{definition}{Definition}[section]
\newtheorem{remark}{Remark}[section]

\numberwithin{equation}{section}
\numberwithin{figure}{section}
\numberwithin{table}{section}



\journal{Journal of Computational and Applied Mathematics}

\begin{document}
\begin{frontmatter}
\title{Clustering-based Low Rank Approximation Method}

\author{Yujun Zhu\fnref{addr0}}
\ead{yujun_zhu@hust.edu.cn}

\author{Jie Zhu\fnref{addr0}}
\ead{zj2021@hust.edu.cn}

\author{Hizba Arshad \fnref{addr0}}
\ead{hizbaarshad96@gmail.com}

\author{Zhongming Wang\fnref{addr2}}
\ead{zwang6@fiu.edu}

\author{Ju Ming\fnref{addr0,addr1}\corref{cor1}} 
\ead{jming@hust.edu.cn}
\cortext[cor1]{Corresponding author}

\address[addr0]{School of Mathematics and Statistics, Huazhong University of Science and Technology, Wuhan, China}

\address[addr1]{Hubei Key Laboratory of Engineering Modeling and Scientific Computing, Huazhong University of Science and Technology, Wuhan, China}

\address[addr2]{Department of Mathematics and Statistics, Florida International University, Miami, FL, USA}

\begin{abstract}
{
We propose a clustering-based generalized low rank approximation method, which takes advantage of appealing features from both the generalized low rank approximation of matrices (GLRAM) and cluster analysis. It exploits a more general form of clustering generators and similarity metrics so that it is more suitable for matrix-structured data relative to conventional partitioning methods. In our approach, we first pre-classify the initial matrix collection into several small subset clusters and then sequentially compress the matrices within the clusters. This strategy enhances the numerical precision of the low-rank approximation. In essence, we combine the ideas of GLRAM and clustering into a hybrid algorithm for dimensionality reduction. The proposed algorithm can be viewed as the generalization of both techniques. Theoretical analysis and numerical experiments are established to validate the feasibility and effectiveness of the proposed algorithm.
}
\end{abstract}


\end{frontmatter}

\section{Introduction}

In the era of rapid technological advancement, the development of data collection tools, techniques, and storage capabilities has led to the generation of massive volumes of data. Modern high dimensional data (HDD) is continuously evolving in various formats including text, digital images, speech signals and videos. However, it has become a challenging task to analyze and process HDD due to its expensive computational complexity. As the dimensionality of the data increases, the efficiency and effectiveness of algorithms tend to deteriorate exponentially, a phenomenon commonly referred to as the \textit{curse of dimensionality}. This results in higher computational costs, longer processing times, and greater resource demands. Dimensionality reduction techniques provide an efficient solution by reducing the number of dimensions and extracting relevant features for analysis, thereby alleviating the computational burden.

The problem of dimensionality reduction has recently received significant attention in areas such as face recognition \cite{chen2012low, huang2012nonlinear, mandal2009curvelet, marjan2021pca, turk1991eigenfaces, zhao2003face}, machine learning \cite{castelli2003csvd, chia2022sampling, entezari2020all, indyk2019learning, srebro2003weighted} and information retrieval \cite{berry1995using, he2023similarity, jin2008ranking, liu2009learning, yu2020active}. Its primary goal is to obtain a more compact data representation with minimal loss of information. During the dimensionality reduction process, the original features are projected into a lower-dimensional space, effectively reducing the redundancy of features and mitigating the curse of dimensionality in high-dimensional data. There are numerous essential dimensionality reduction techniques based on vector selection, i.e., the vector space model \cite{turk1991eigenfaces, zhao2003face}. It treats each sample as a vector, and the entire dataset is represented as a matrix, with the columns corresponding to the data points. One of the most well-known techniques based on the vector space model is low rank matrix approximation using singular value decomposition (SVD) \cite{klema1980singular}. An appealing property of SVD is that it can achieve the smallest reconstruction error, given the same rank, in terms of Euclidean distance \citep{eckart1936approximation}, making it widely applicable in areas such as signal processing \cite{mademlis2018regularized, sahidullah2016local}, big data \cite{bjornsson1997manual, falini2022review,  tanwar2018dimensionality}, etc.

The standard SVD-based low-rank approximation of matrices operates on a single data matrix at one time, which limits its application in vector-format data. Additionally, SVD faces high computational and space complexities, particularly when dealing with high-dimensional or tensor-format data. Many attempts have been made in order to overcome these obstacles. Ye proposed the generalized low rank approximation (GLRAM) method \cite{ye2004generalized}, which is proved to have less CPU elapsing time than traditional SVD. It transfroms the original dataset into a smaller subspace by employing a pair of shared left- and right- projectors. Numerous improvements and extensions of the GLRAM method have also been developed. Liu provided a non-iterative solution of the corresponding minimization program \cite{liu2006non}. A simplified form was present in \cite{lu2008simplified}, which could further lessen the GLRAM computational complexity. To deal with large sparse noise or outliers, the robustness of such low rank approximation approach was discussed in \cite{nakouri2016robust, shi2015robust, zhao2016robust}. The combinations of GLRAM, machine learning and deep neutral techniques were also applied in pattern recognition and image degradation \cite{ahmadi2020generalized, indyk2021few, ren2018deep}. We refer to \cite{zhu2024low} for the application of low rank approximation in solving stochastic partial differential equations and optimal control problems.

Despite GLRAM's capability to preserve spatial relationships and relatively low computational complexity, its accuracy largely depends on the intrinsic correlations within the original dataset. More precisely, GLRAM is restricted to the search space with a small sample size or high correlations \cite{liu2010generalized, ye2004generalized}. In this article, we aim to broaden the universality and applicability of the GLRAM method so that it can perform well in real-world scenarios, where large datasets and complex correlation structures between data matrices are common. To achieve this goal, we integrate clustering techniques with the existing GLRAM method. 

Data clustering is a form of unsupervised classification to analyze similarities and dissimilarities of intrinsic characteristics among different samples. It has been widely adopted in fields such as statistical analysis \cite{dalmaijer2022statistical, srivastava2005statistical, von2005towards}, pattern recognition \cite{ali2017k,schwenker2014pattern}, and learning theory \cite{ezugwu2022comprehensive, petegrosso2020machine}. During the clustering process, the original dataset is divided into several clusters, with samples exhibiting greater similarity grouped together. K-means clustering method \cite{lloyd1982kmeans, macqueen1967some, selim1984k} is one of the most popular centroid-based algorithms, where samples are reassigned based on their distances to the cluster centroids. K-means has appealing properties such as ease of implementation, fast computation and suitability for high-dimensional data. In computational science, it is also known as a special case of centroidal Voronoi tessellations (CVT), which generalizes the method to continuous data space \cite{du2006convergence, du1999centroidal, du2002centroidal, du2010advances, du2005optimal, du2002numerical}.

In this article, we propose a clustering-based generalized low rank approximation method, where the GLRAM and K-means clustering techniques are combined to define a generalization of GLRAM. The novel dimensionality reduction technique leverages the strengths of both GLRAM and K-means. On one hand, it adopts a more general form of clustering centroid and distance, allowing it to be applied to data with more complex formats. On the other hand, our algorithm pre-classifies the original dataset and separately compresses sampling matrices by cluster, thereby improving the accuracy and efficiency of low rank approximation.

The rest of the article is organized as follows. Section 2 respectively reviews two low rank approximation methods, SVD and GLRAM, and some of their properties. In Section 3, we briefly introduce the concept of K-means clustering and its algorithmic process. Section 4 presents our algorithm of clustering-based generalized low rank approximation of matrices, compares it with the standard SVD and GLRAM methods, and provides some theoretical analysis. In Section 5, we illustrate our numerical results, which validate the feasibility and the effectiveness of our proposed algorithm. Finally, Section 6 provides some concluding remarks.

\section{Low rank approximation of matrices}

The nature of low rank approximation of matrices is dimensionality reduction. It aims to obtain a more compact data representation with limited loss of information, so that solving a given problem on the lower-rank alternative approximates the solution on original matrices, which leads to less computational and space complexity. 

Mathematically, the optimal rank-$k$ approximation ($k \le N$) of a matrix $\mathbb{A} \in \mathbf{R}^{N \times N}$ under the Frobenius norm is formulated as a rank-constrained minimization problem (\ref{eq2.1}):  find a matrix $\widetilde{\mathbb{A}} \in \mathbf{R}^{N \times N}$ such that

\begin{equation} \label{eq2.1}
	\widetilde{\mathbb{A}}^* \: = \: \mathop{\arg\min}\limits_{rank(\widetilde{\mathbb{A}}) \; = \; k} \enspace \Vert \mathbb{A} - \widetilde{\mathbb{A}} \Vert_F,
\end{equation}

\noindent while $\Vert \mathbb{A} - \widetilde{\mathbb{A}} \Vert_F$ is called as the \textit{reconstruction error} of the approximation.

\subsection{Singular value decomposition}

Singular value decomposition is one of the most well-known techniques in low rank approximation of matrices. It is the factorization of $\mathbb{A}$ into the product of three matrices $\mathbb{A} = \mathbb{U} \Sigma \mathbb{V}^T \in \mathbf{R}^{N \times N}$, where the columns of $\mathbb{U}, \mathbb{V}$ are orthonormal and $\Sigma$ is a diagonal matrix with positive real entries. Such decomposition always exists for any complex matrix. 

An appealing property of low rank approximation of matrices via SVD is that it can achieve the smallest reconstruction error among all approximations with the constraint of the same rank in Euclidean distance \cite{eckart1936approximation}. In other words, the minimization program in Eq. (\ref{eq2.1}) admits an analytical solution in terms of its truncated singular value decomposition (TSVD), as stated in the following theorem. 

\begin{theorem}[Eckart-Young Theorem \cite{eckart1936approximation}] \label{th2.1} 
	
Let $\mathbb{A} = \mathbb{U} \Sigma \mathbb{V}^T \in \mathbf{R}^{N \times N}$ be the SVD of $\mathbb{A}$, and let $\mathbb{U}, \Sigma$ and $\mathbb{V}$ partitioned as follows:
	
\begin{equation} 
	\mathbb{U} \: := \: \begin{bmatrix} \mathbb{U}_1 & \mathbb{U}_2 \end{bmatrix}, \quad \Sigma \: := \: \begin{bmatrix} \Sigma_1 & 0 \\ 0 & \Sigma_2 \\ \end{bmatrix}, \quad and \quad \mathbb{V} \: := \: \begin{bmatrix} \mathbb{V}_1 & \mathbb{V}_2  \end{bmatrix},
	\nonumber
\end{equation}
	
\noindent where $\mathbb{U}_1, \mathbb{V}_1 \in \mathbf{R}^{N \times k}$ and $\Sigma_1 \in \mathbf{R}^{k \times k}$. Then the rank-k matrix, obtained from the TSVD, $\widetilde{\mathbb{A}^*} = \mathbb{U}_1 \Sigma_1 \mathbb{V}_1^T$, satisfies that
	
\begin{equation} \label{eq2.2}
	\Vert \mathbb{A} - \widetilde{\mathbb{A}}^* \Vert_F \: = \: \mathop{\min}\limits_{rank(\widetilde{\mathbb{A}}) \leq k} \: \Vert \mathbb{A} - \widetilde{\mathbb{A}} \Vert_F.
\end{equation}
	
\noindent The minimizer $\widetilde{\mathbb{A}}^*$ is unique if and only if $\sigma_{k+1} \neq \sigma_k$.
	
\end{theorem}

The main drawback of SVD is that the application of the technique in large matrices encounters practical limits both in time and space aspects, due to the expensive computation. When applied to high-dimensional data, such as images and videos, we will quickly run up against practical computational limits. Therefore, the method of generalized low rank approximations of matrices \cite{ye2004generalized} is proposed to alleviate the high SVD computational cost.

\subsection{Generalized low rank approximation of matrices}

As the generalization of SVD, the GLRAM algorithm still uses two-dimensional matrix patterns and it simultaneously carries out tri-factorizations on a collection of matrices. Let $\{\mathbb{A}_i\}_{i=1}^N \in \mathbf{R}^{r \times c}$ be $N$ data matrices. The GLRAM goal is to seek orthonormal matrices $\mathbb{L} \in \mathbf{R}^{r \times k_1}$ and $\mathbb{R} \in \mathbf{R}^{k_2 \times c}$ with $k_1 \ll r$ and $k_2 \ll c$, such that 

\begin{equation} \label{eq2.3}
\min_{\substack{
			\mathbb{L} \in \mathbf{R}^{r \times k_1}, \mathbb{L}^T \mathbb{L} = I_{k_1}  \\
			\mathbb{R} \in \mathbf{R}^{c \times k_2}, \mathbb{R}^T \mathbb{R} = I_{k_2} \\
		    \mathbb{M}_i \in \mathbf{R}^{k_1 \times k_2}, i = 1,2,...,N }}
\quad \sum_{i=1}^N \Vert \mathbb{A}_i - \mathbb{L} \mathbb{M}_i \mathbb{R}^T \Vert_F^2,
\end{equation}

\noindent where it holds that $\mathbb{M}_i $ is the compressed version of $\mathbb{A}_i$,

\begin{equation} \label{eq2.4}
\mathbb{M}_i = \mathbb{L}^T \mathbb{A}_i \mathbb{R}, \quad i = 1,\cdots N.
\end{equation}

Here $\tilde{\mathbb{A}}_i = \mathbb{L} \mathbb{M}_i \mathbb{R}^T$ is called a reconstructed matrix of $\mathbb{A}_i$, $i=1,2, ... ,N$. In order to make the above reconstruction error as small as possible, Ye suggested that the values of $k_1$ and $k_2$ should be equal (denoted by $k$) \cite{ye2004generalized}. Although there is no strict theoretical proof for this, Liu et al. \cite{liu2010generalized}  provided some theoretical support for such a choice.

Furthermore, it was shown that the above minimization problem is mathematically equivalent to the following optimization problem

\begin{equation} \label{eq2.5}
	\max_{\substack{
			\mathbb{L} \in \mathbf{R}^{r \times k}, \mathbb{L}^T \mathbb{L} = I_{k}  \\
			\mathbb{R} \in \mathbf{R}^{k \times c}, \mathbb{R}^T \mathbb{R} = I_{k} }}
	\quad \sum_{i=1}^N \Vert  \mathbb{L}^T \mathbb{A}_i \mathbb{R} \Vert_F^2.
\end{equation}

Unfortunately, both Eq. (\ref{eq2.3}) and Eq. (\ref{eq2.5})  have no closed-form solutions in general, and we have to solve them by using some iterative algorithms, as the following theorem indicates.

\begin{theorem}[\cite{ye2004generalized}] \label{th2.2} 
	
Let $\mathbb{L}$ and $\mathbb{R}$ be the optimal solutions to the problem in Eq. (\ref{eq2.3}), then that holds

(i) Given matrix $\mathbb{R}$, $\mathbb{L}$ constitutes the $k$ eigenvectors of the matrix corresponding to the largest $k$ eigenvalues:

\begin{equation} \label{eq2.6}
\mathbb{N}_L(\mathbb{R}) = \sum_{i=1}^N \mathbb{A}_i \mathbb{R} \mathbb{R}^T \mathbb{A}_i^T.
\end{equation}

(ii) Given matrix $\mathbb{L}$, $\mathbb{R}$ constitutes the $k$ eigenvectors of the matrix corresponding to the largest $k$ eigenvalues:

\begin{equation}\label{eq2.7}
\mathbb{N}_R(\mathbb{L}) = \sum_{i=1}^N \mathbb{A}_i^T \mathbb{L} \mathbb{L}^T \mathbb{A}_i.
\end{equation}
	
\end{theorem}

That is, we achieve $\mathbb{L}$ and $\mathbb{R}$ by alternately solving $k$ dominant eigenvectors of $\mathbb{N}_L$ and $\mathbb{N}_R$, respectively. This yields an iterative algorithm for solving the projection matrices summarized in Algorithm \ref{alg1}.

\begin{algorithm}[h]
	\caption{The generalized low rank approximations of matrices algorithm} \label{alg1}
	\begin{algorithmic}[1]
		
		\Require
		Matrices $\{\mathbb{A}_i\}_{i=1}^M$, and the dimension $k$
		\Ensure
		Matrices $\mathbb{L}, \mathbb{R}$ and $\{\mathbb{M}_i\}_{i=1}^M$
		
		\State Obtain initial $\mathbb{L}_0$ and set iteration count $t = 1$.
		
		\While{not convergent}
		
		\State Form matrix $\mathbb{N}_R$ by Eq. (\ref{eq2.6}).
		
		\State Compute the $k$ eigenvectors $\{\phi_R^j\}_{j=1}^k$ of $\mathbb{N}_R$ corresponding to the largest $k$ eigenvalues, and let $\mathbb{R}_{t} = \left[ \phi_R^1,...,\phi_R^k \right]$.
		
		\State Form matrix $\mathbb{N}_L$ by Eq. (\ref{eq2.7}).
		
		\State Compute the $k$ eigenvectors $\{\phi_L^j\}_{j=1}^k$ of $\mathbb{N}_L$ corresponding to the largest $k$ eigenvalues, and let $\mathbb{L}_{t} = \left[ \phi_L^1,...,\phi_L^k \right]$.
		
		\EndWhile
		
		\State \Return Matrices $\mathbb{L} = \mathbb{L}_{t}, \mathbb{R} = \mathbb{R}_{t}, \mathbb{M}_i = \mathbb{L}^T \mathbb{A}_i \mathbb{R}, i = 1,2,...,N$.
		
	\end{algorithmic}
\end{algorithm}

Theorem \ref{th2.2} also implies that the matrix updates in Lines 4 and 6 in Algorithm \ref{alg1} do not decrease the value of  $\sum_{i=1}^N \Vert  \mathbb{L}^T \mathbb{A}_i \mathbb{R} \Vert_F^2$, since the computed $\mathbb{L}$ and $\mathbb{R}$ are locally optimal. Hence the value of the Root Mean Square Reconstruction Error (RMSRE),

\begin{equation} \label{eq2.8}
	RMSRE(M) \: := \: \sqrt{\frac{1}{M} \sum_{m=1}^M \Vert \mathbb{A}_m - \mathbb{L} \mathbb{M}_m \mathbb{R}^T \Vert_F^2},
\end{equation}

\noindent does not increase. The convergence of Algorithm \ref{alg1} follows as stated in the following theorem, since RMSRE is bounded from below by 0. This convergence property is summarized in the theorem below.

\begin{theorem}[\cite{ye2004generalized}] \label{th2.3} 
	
Algorithm \ref{alg1} monotonically non-increases the value of the Root Mean Square Reconstruction Error (RMSRE) as defined in Eq. (\ref{eq2.8}), and therefore it converges in the limit.

\end{theorem}

\begin{figure}[ht]
\centering
\includegraphics[scale=0.36]{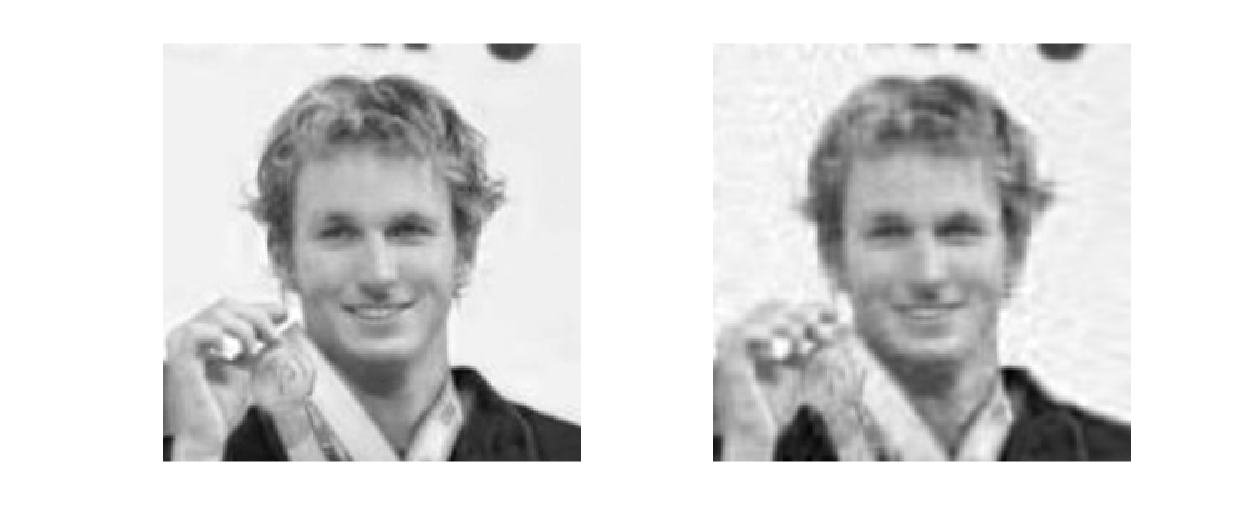}
\caption{Raw image from LFW dataset (left) and image compressed by GLRAM with reduction ratio $10\%$ (right).}
\label{fig2.1}
\end{figure}

GLRAM treats data as native two-dimensional matrix patterns rather than conventionally using vector space. It has also been verified experimentally to have better compression performance and enjoy less computational complexities than traditional dimensionality reduction techniques \cite{ye2004generalized}. Due to these appealing properties, GLRAM has naturally become a good choice in image compression and retrieval, where each image is represented in its native matrix form. Fig \ref{fig2.1} reveals the nice compression performance of GLRAM, which could well capture the main features of the face image while requiring small storage. Here the image is randomly chosen from the LFW dataset \cite{huang:inria-00321923}.

\section{Clustering Technique}

Clustering is a tool to analyze similarities and dissimilarities between different samples, which partitions the entire sample collection into several clusters according to the similarities between samples. K-means clustering \cite{macqueen1967some} is one of the simplest and most popular centroid-based clustering algorithms, which aims to reassign a set of data points according to their distance with the centroid of each cluster.

Consider a set of observations $\{x_1, x_2, ..., x_N\} \in  \mathbf{R}^{d}$, where each observation is a $d$-dimensional real vector, K-means method aims to partition these $N$ observations into $K$ clusters $\mathbf{V} = \{\mathbf{V}_1, \mathbf{V}_2, ..., \mathbf{V}_K\}$ with $K \le N$, so as to minimize the within-cluster sum of square error (SSE) (i.e. within-cluster variance). Mathematically, the objective function is stated as below:

\begin{equation} \label{eq3.1}
	\mathop{\arg\min}_{\mathbf{V}} \sum_{j=1}^K \sum_{x \in \mathbf{V}_j} \Vert x - \mu_j \Vert^2 = \mathop{\arg\min}_{\mathbf{V}} \sum_{j=1}^K \left\vert\mathbf{V}_j\right\vert Var(\mathbf{V}_j).
\end{equation}

\noindent Here $\mu_j$ is the arithmetic mean of points in $\mathbf{V}_j$ and it is also called the centroid of the convex set $\mathbf{V}_j$, i.e.,

\begin{equation} \label{eq3.2}
	\mu_j = \frac{1}{\left\vert \mathbf{V}_j \right\vert} \sum_{x \in \mathbf{V}_j} x,
\end{equation} 

\noindent where $\left\vert \mathbf{V}_j \right\vert$ is the size of $\mathbf{V}_j$, and $\Vert \cdot \Vert$ is the usual Euclidean distance (the $L^2$ metric) in $\mathbf{R}^d$. 

The following theorem shows that for any fixed clustering, the sum of square error is minimized when the centroid associated with each cluster is the arithmetic mean of the set of points assigned to that cluster, which validates the effectiveness of the updating formulation in Eq.(\ref{eq3.2}) and also the convergence of the Lloyd's K-means method.

\begin{theorem}[\cite{burkardt2009k}] \label{th3.1} 
	
Given a set of data $\mathbf{V}$ with the size of  $\vert \mathbf{V}\vert$ and for each datapoint $x \in \mathbf{R}^d$. Let $x$ be an arbitrary point in the same d-dimensional space. Then it holds
	
\begin{equation} \label{eq3.3}
	\sum_{x_i \in \mathbf{V}} \Vert x_i - x \Vert^2 \ge \sum_{x_i \in \mathbf{V}} \Vert x_i - \mu \Vert^2,
\end{equation} 

\noindent where $\mu$ the arithmetic mean of these points

$$\mu = \frac{1}{\vert \mathbf{V}\vert} \sum_{x \in \mathbf{V}} x.$$
	
\end{theorem}

\begin{remark}

K-means is all about the analysis-of-variance paradigm. Both univariate and multivariate analysis of Variance (ANOVA) are based on the fact that the sum of squared deviations about the cluster centroid is comprised of such scatter about the group centroids and the scatter of those centroids about the whole cluster.

\begin{theorem}[\cite{ingrassia2020cluster}] \label{th3.2} 
	
	Given a set of data $\mathbf{V}$ with the size of  $\vert \mathbf{V}\vert$ and the number of clusters $K$, the within-cluster sum of squares total (SST) corresponding to the K-means method can be decomposed into the SSE and between-clusters sum of square (SSB) described as below, which is also called between-clusters variance, 
	
	\begin{equation} \label{eq3.4}
		SST = SSE + SSB,
	\end{equation} 
	
	\noindent where
	
	\begin{equation} \label{eq3.5}
		\begin{aligned}
			& SST = \sum_{j=1}^K \sum_{x \in \mathbf{V}_j} \Vert x - \mu \Vert^2,\\
			& SSB = \sum_{j=1}^K \vert \mathbf{V}_j \vert \Vert \mu_j - \mu \Vert^2,
		\end{aligned}
	\end{equation} 
	
	\noindent where $\mu = \frac{1}{\vert \mathbf{V}\vert} \sum_{i=1}^N x_i$ is the centroid of $\mathbf{V}$.
	
\end{theorem}

Therefore, minimizing the objective function in Eq. (\ref{eq3.1}) will necessarily maximize the distance between clusters.
	
\end{remark}

\begin{remark}
The k-means clustering technique can be seen as partitioning the space into Voronoi cells defined in \cite{okabe2009spatial},

\begin{equation} \label{eq3.6}
	\mathbf{V}_j := \left\lbrace  x \in \mathbf{V} \vert \Vert x - \mu_j \Vert \le \Vert x - \mu_i \Vert, \quad \forall i = 1,2,...,K \quad and \quad i \neq j \right\rbrace ,
\end{equation} 

\noindent where $\mathbf{V} = \cup_{j=1}^K \mathbf{V}_j$ and $\mathbf{V}_i \cap \mathbf{V}_j = \emptyset$ for any $i \neq j$, so that $\{\mathbf{V}_j\}_{j=1}^K$ is the Voronoi tessellation of $\mathbf{V}$ associated with the centroid set $\{\mu_j\}_{j=1}^K$.

For each pair of two centroids, there is a line that connects them. Perpendicular to this line, there is a line, plane or hyperplane (depending on the dimensionality) that passes through the middle point of the connecting line and divides the space into two separate subspaces, as depicted in Figure \ref{fig3.1}. The k-means clustering therefore partitions the space into $K$ subspaces for which $\mu_i$ is the nearest centroid for all included elements of the subspace \cite{faber1994clustering}.

\begin{figure}[ht]
	\centering
	\includegraphics[scale=0.3]{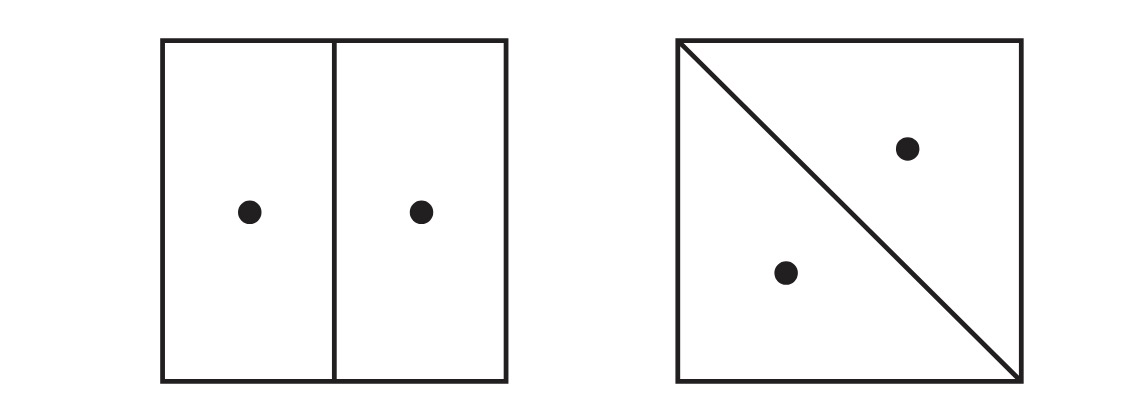}
	\caption{Two centroidal Voronoi tessellations of a square. The points $\mu_1$  and $\mu_2$ are the centroids of the rectangles on the left or of the triangles on the right.}
	\label{fig3.1}
\end{figure}
\end{remark}

There are several algorithms known for determining K-means clusters of a given dataset in \cite{du1999centroidal, faber1994clustering, macqueen1967some, lloyd1982kmeans}. One of famous approaches is a deterministic algorithm known as Lloyd’s method \cite{lloyd1982kmeans}, in which $K$ data points are initially chosen as centroids at random, and all samples are then assigned to their closest centroid following the recalculation of centroids as arithmetic mean over all their assigned groups. Lloyd's algorithm loops over Lines 2-3 until convergency, which is summarized in Algorithm \ref{alg2}. There are also various other methods based on the minimization properties of K-means clustering. 

\begin{algorithm}[h]  
	\caption{Lloyd's K-means clustering algorithm}
	\begin{algorithmic}[1]
		
		\State Randomly select $N$ centroids.
		
		\State Measure the similarities between all the data matrices and the centroids and assign them to different clusters.
		
		\State Update the centroids in each cluster.
		
		\State Repeat Steps 2 - 3 until convergence or maximal iteration count is reached.
		
	\end{algorithmic} \label{alg2}
\end{algorithm}

K-means clustering is a centroid-based technique, and it is essential to figure out how to obtain centroids and measure the similarities. Traditionally, we pick out several sample matrices as the initial centroids and replace them with more suitable samples according to their Euclidean distance. However, the raw data generally have large dimensions and it is computationally expensive to acquire their norms. Therefore, we aim to find a generalized form to measure the similarities of matrices and update the centroids in each cluster.

\section{Generalized low rank approximation of matrices based on clustering methods}

In this section, we would like to combine GLRAM and the K-means clustering methods, in order to take advantage of appealing properties of both techniques. Here we need to first extend the standard K-means clustering settings, i.e., the more general notions of centroid and distance, which are two main innovations in our clustering procedure.

First, given a sequence of matrices $\mathbf{V}:=\{\mathbb{A}^1, \cdots, \mathbb{A}^n\} \in \mathbf{R}^{r \times c}$, we define the \textit{rank-$k$ generalized centroid} ($k \le \min\{r,c\}$) of such observatory set is two lower-dimensional matrices $\mathbb{L} \in \mathbf{R}^{r \times k}, \mathbb{R}\in \mathbf{R}^{c \times k}$ with orthonormal columns, subject to

\begin{equation}\label{eq4.1}
	\begin{aligned}
		\min \: \sum_{\mathbb{A}^i \in \mathbf{V}} \:   \Vert \mathbb{A}^i \;-\;   \mathcal{P}(\mathbb{A}^i, \mathbf{V}) \Vert_F^2 \:=\: \min_{\substack{
				\mathbb{L} \in \mathbf{R}^{r \times k}, \mathbb{L}^T \mathbb{L} = I_{k}  \\
				\mathbb{R} \in \mathbf{R}^{c \times k}, \mathbb{R}^T \mathbb{R} = I_{k} \\
				\mathbb{M}^i \in \mathbf{R}^{k \times k}, i = 1,2,...,n }}\:\sum_{i=1}^n \:  \Vert \mathbb{A}^i \;-\;  \mathbb{L} \mathbb{M}^i \mathbb{R}^T \Vert_F^2,
	\end{aligned}
\end{equation}

\noindent which also corresponds to the objective function of the generalized low rank approximation of the set $\mathbf{V}$ in Eq. (\ref{eq2.3}). For notational simplicity, the centroid is denoted as the generalized left- and right- low rank approximation matrix pairs $(\mathbb{L}, \mathbb{R})$, which could be obtained by Algorithm \ref{alg1}. $(\mathbb{L}, \mathbb{R})$ act as the two-sided linear transformations on the data in matrix form, and $\mathcal{P}$ could be seen as the projection operator into the rank-$k$ subspace of $\mathbf{V}$, which also connects with our new centroids $(\mathbb{L}, \mathbb{R})$. We do not employ the arithmetic mean in Eq. (\ref{eq3.2}) as the centroid anymore. However, our approach still has similarities with K-means clustering in the manner of determining centroids. We find optimal rank-$k$ projections of the initial observations related to $(\mathbb{L}, \mathbb{R})$, and equally, the centroids $(\mathbb{L}, \mathbb{R})$ also contain information about the samples within the matrix set.  

Second, we manipulate different variable spaces and similarity metrics. K-means clustering faces observations in vector forms, and its square of the similarity metric between two vectors $\bm{x}, \bm{y} \in \mathbf{R}^d$ is defined by

\begin{equation}\label{eq4.2}
\begin{aligned}
	\delta^2 (\bm{x}, \bm{y}) \: & = \: \Vert \bm{x} \; - \; \bm{y} \Vert_2^2 \\
	 & = \: \sum_{l=1}^d\; (x_l \; - \; y_l)^2,
\end{aligned}
\end{equation}

\noindent which is also called the Euclidean distance between vectors $\bm{x}$ and $\bm{y}$.  To generalize the K-means approach, the observatory sequence of our new approach is essentially concerned with two-dimensional matrices $\mathbb{A}, \mathbb{B} \in \mathbf{R}^{r \times c}$ and we adopt the Frobenius norm given as follows to determine the reconstruction error.

\begin{equation}\label{eq4.3}
\begin{aligned}
    \delta^2 (\mathbb{A}, \mathbb{B}) \:& = \: \Vert \mathbb{A}  \; - \;  \mathbb{B} \Vert_F^2 \\
	& = \: \sum_{i=1}^r \sum_{j=1}^c  \; (a_{ij}  \; - \;  b_{ij})^2.
\end{aligned}
\end{equation}

In this setting, two-dimensional matrix data is equivalent to one-dimensional vector data. According to the distance frameworks in Eq. (\ref{eq4.2}) and (\ref{eq4.3}),  we can convert the matrices $\mathbb{A}, \mathbb{B}$ to the $d$-dimensional ($d=r \cdot c$) vectors $\bm{a} = vec(\mathbb{A})$ and $\bm{b} = vec(\mathbb{B})$, where we sequentially concatenate the columns of $\mathbb{A}, \mathbb{B}$. Meanwhile, the relationship between the two similarity metrics is directly derived by the definitions of the Frobenius norm and the Euclidean norm that: the Frobenius norm of the $r \times c$ difference matrix $\mathbb{A} - \mathbb{B}$ is identical to the Euclidean distance of the $d$ vectors $\bm{a}$ with $d = r \cdot c$.
	
\begin{equation}
\begin{aligned}
	\delta (\mathbb{A}, \mathbb{B})  \:&=\: \sum_{i=1}^r \sum_{j=1}^c  \; (a_{ij}  \; - \;  b_{ij})^2 \\
	&=\: \sum_{l=1}^d \; (a_l \; - \; b_l)^2 \enspace=\: \delta (\bm{a}, \bm{b}).
\nonumber
\end{aligned}
\end{equation}

More generally, the square of the distance between a matrix $\mathbb{A}$ and the centroid $(\mathbb{L}, \mathbb{R})$ of the matrix set $\mathbf{V}$ is defined by

\begin{equation}\label{eq4.4}
	\delta^2 \left( \mathbb{A}, (\mathbb{L}, \mathbb{R})\right)   \quad:=\: \delta^2 \left( \mathbb{A},  \mathcal{P}(\mathbb{A}^i, \mathbf{V})\right)
		\:=\: \Vert \mathbb{A} \;-\;  \mathbb{L} \mathbb{M} \mathbb{R}^T \Vert_F^2,
\end{equation}

\noindent where the middle low-dimensional matrix is determined by $\mathbb{M} = \mathbb{L}^T \mathbb{A} \mathbb{R}$. Then, given a matrix sequence $\mathbf{V} = \{\mathbb{A}^i\}$ and a set of rank-$k$ generating matrices $\{(\mathbb{L}_j,\mathbb{R}_j)\}_{j=1}^K$, the \textit{rank-$k$ generalized clusters} of $\mathbf{V}$ is given by
	
\begin{equation}\label{eq4.5}
\begin{aligned}
	\mathbf{V_j} \:=\: \{\mathbb{A}^i \in \mathbf{V} \;\vert\;  \delta^2\left( \mathbb{A}^i, (\mathbb{L}_j, \mathbb{R}_j)\right)  \;\le\;  \delta^2\left( \mathbb{A}^i, (\mathbb{L}_l, \mathbb{R}_l)\right) , \: \forall \;l \neq j\} \quad \forall \;j = 1,\cdots,K,
\end{aligned}
\end{equation}

\noindent where $\mathbf{V_j}$ is the $j$-th subset divided by the original dataset $\mathbf{V}$, and it satisfies that $\mathbf{V} = \cup_{j=1}^K \mathbf{V}_j$ and $\mathbf{V}_i \cap \mathbf{V}_j = \emptyset$ for any $i \neq j$.

Now, we are ready to define our clustering-based generalized low rank approximation of matrices (CGLRAM) method. Given a sequence of feature matrices $\mathbf{V} := \{\mathbb{A}^i\}_{i=1}^N$ with each matrix $\mathbb{A}^i \in \mathbf{R}^{r \times c}$, we aim to find our rank-$k$ generalized clusters $\{\mathbf{V_j}\}_{j=1}^K$ satisfying Eq. (\ref{eq4.5}) and rank-$k$ generalized centroids $\{(\mathbb{L}_j,\mathbb{R}_j)\}_{j=1}^K$ given by Eq. (\ref{eq4.1}). Actually, the rank-$k$ generalized clusters and centroids are in exact correspondence. In other words, the tessellation $\{\mathbf{V_j}\}_{j=1}^K$ is called a rank-$k$ generalized clusters if and only if their corresponding $\{(\mathbb{L}_j,\mathbb{R}_j)\}_{j=1}^K$ are themselves the rank-$k$ generalized centroids of the $\mathbf{V}_j$'s, as stated as the following theorem.

\begin{theorem}\label{th4.1}
	
Let $\mathbb{L}_j, \mathbb{R}_j$ be the optimal solution to the following minimization problem
	
\begin{equation} \label{eq4.6}
	RE(\mathbf{V}_j) \;  := \; \min \;\sum_{\mathbb{A}^i \in \mathbf{V}_j} \Vert  \mathbb{A}^i - \mathbb{L}_j \mathbb{M}^i \mathbb{R}_j^T \Vert_F^2,
\end{equation}
	
\noindent then $\{\mathbb{L}_j\}_{j=1}^K, \{\mathbb{R}_j\}_{j=1}^K$ are the rank-$k$ generalized centroids of the whole set $\mathbf{V}$, and it holds that
	
\begin{equation} \label{eq4.7}
	\mathbb{M}^i \; = \; \mathbb{L}^T \mathbb{A}^i \mathbb{R}, \quad \forall \mathbb{A}^i \in \mathbf{V}_j.
\end{equation} 
	
\end{theorem}

\begin{proof}
	
By the definition in Eq. (\ref{eq4.6}), $RE(\cdot ): \mathbf{V}_j \rightarrow \mathbf{R}$ is a mapping from the data matrices set to the real space. Suppose $\mathbb{L}_j, \mathbb{R}_j$ are the minimizers of subsets $\mathbf{V}_j $, for $j = 1,...,N$, i.e. let $\widetilde{\mathbb{L}}_j, \widetilde{\mathbb{R}}_j $ be arbitrary matrices with orthonormal columns satisfying that $\widetilde{\mathbb{L}}_j \in \mathbf{R}^{r \times k}$ and $\widetilde{\mathbb{R}}_j \in \mathbf{R}^{c \times k}$, it must hold that
	
$$ \sum_{\mathbb{A}^i \in \mathbf{V}_j} \Vert  \mathbb{A}^i - \mathbb{L}_j \mathbb{M}^i \mathbb{R}_j^T \Vert_F^2 \;  \le \; \sum_{\mathbb{A}^i \in \mathbf{V}_j} \Vert  \mathbb{A}^i - \widetilde{\mathbb{L}}_j \mathbb{M}^i \widetilde{\mathbb{R}}_j^T \Vert_F^2, \quad j = 1,2,...,K.$$
	
Therefore, by summing the inequalities above together, we have 
	
$$ \sum_{j=1}^K\sum_{\mathbb{A}^i \in \mathbf{V}_j} \Vert  \mathbb{A}^i - \mathbb{L}_j \mathbb{M}^i \mathbb{R}_j^T \Vert_F^2 \;  \le \; \sum_{j=1}^K \sum_{\mathbb{A}^i \in \mathbf{V}_j} \Vert  \mathbb{A}^i - \widetilde{\mathbb{L}}_j \mathbb{M}^i \widetilde{\mathbb{R}}_j^T \Vert_F^2,$$
	
\noindent where the equality in the formulation above is achieved if and only if $\widetilde{\mathbb{L}}_j = \mathbb{L}_j , \widetilde{\mathbb{R}}_j = \mathbb{R}_j$, for $j = 1,...,K$. Therefore, $\{\mathbb{L}_j\}_{j=1}^K, \{\mathbb{R}_j\}_{j=1}^K$ are the rank-$k$ generalized centroids of the set $\mathbf{V}$. By Theorem 3.1 in \cite{ye2004generalized}, the key step for the minimization in Eq. (\ref{eq4.6}) is the computation of the orthonormal transformations $\mathbb{L}_j$ and $\mathbb{R}_j$, and $\mathbb{M}^i$ is uniquely determined by $\mathbb{L}_j$ and $\mathbb{R}_j$ with Eq. (\ref{eq4.7}). Hence we derive the original formula.
	
\end{proof}

To recapitulate, CGLRAM can be viewed as a generalization of the standard K-means clustering method, which is applicable to matrix data. The initial matrix sequence $\mathbf{V}$ is divided into $K$ clusters $\{\mathbf{V_j}\}_{j=1}^K$ and the corresponding generators are the generalized left- and right- low rank matrices $\{(\mathbb{L}_j,\mathbb{R}_j)\}_{j=1}^K$. We employ the Forbenius norm to measure the similarities of intrinsic characteristics between samples and clusters. Meanwhile, CGLRAM can also be viewed as a generalization of GLRAM for which the dataset is first pre-classified and we then separately compress the matrices within each cluster, thereby enhancing the accuracy of low rank approximations. Therefore, if $K=1$, then CGLRAM reduces to the standard GLRAM of Section 2.

\begin{equation} \label{eq4.8}
	 \min_{\{(\mathbb{L}_j,\mathbb{R}_j)\}_{j=1}^K} \quad \delta^2\left( \mathbb{A}^i, (\mathbb{L}_j, \mathbb{R}_j)\right)  \quad := \enspace
	\min_{\substack{
			\mathbb{L}_j \in \mathbf{R}^{r \times k}, \mathbb{L}_j^T \mathbb{L}_j = I_{k}\\
			\mathbb{R}_j \in \mathbf{R}^{c \times k}, \mathbb{R}_j^T \mathbb{R}_j = I_{k} \\
			\mathbb{M}^i \in \mathbf{R}^{k \times k}, for \: \mathbb{A}^i \in \mathbf{V}_j \\
			j = 1,2,\cdots,N }}
	\quad\sum_{j=1}^K \sum_{\mathbb{A}^i \in \mathbf{V}_j} \Vert \mathbb{A}^i - \mathbb{L}_j \mathbb{M}^i \mathbb{R}_j^T \Vert_F^2.
\end{equation}

\subsection{Algorithm to construct CGLRAM}

In this subsection, we make a natural extension of Lloyd's method for computing the traditional K-means clusters and centroids to our more general CGLRAM background. Let us begin with a given tessellation $\{\mathbf{V_j}\}_{j=1}^K$. It is actually an iterative procedure that decomposes the whole process of finding the centroids $\{(\mathbb{L}_j,\mathbb{R}_j)\}_{j=1}^K$ into a sequence of low rank approximation analysis within each subset of observatory matrices. More precisely, we first construct the new tessellation by measuring and comparing the similarity between samples and initial centroids and then update the centroids by the corresponding generalized two-sided low rank matrices. The concrete procedure of the clustering-based generalized low rank approximation of matrices is summarized in Algorithm \ref{alg3}.

\begin{algorithm}[h]
\caption{Clustering-based Generalized Low Rank Approximation of Matrices (CGLRAM)} 
\begin{algorithmic}[1]
	\Require Matrices $\{\mathbb{A}_i\}_{i=1}^N$, cluster count $K$, and the compressed dimension $k$.
	\Ensure Low rank matrices $\mathbb{L}_j$ , $\mathbb{R}_j$ and $\{\mathbb{M}^i\}_{i=1}^{\vert \mathbf{V}_j \vert}$, for $j = 1,\cdots, K$ s.t. $\{\mathbf{V}_j\}_{j=1}^K$ is a tessellation of $\mathbf{V}$.
		
	\State Randomly divide $K$ matrix samples as initial clustering centriods and compute their generalized low rank left- and right- orthogonal matrices $\mathbb{L}_j, \mathbb{R}_j$ by Algorithm \ref{alg1}.
		
	\State Classify each matrix sample into the clusters corresponding to the centroid with the smallest distance between sample matrix $\mathbb{A}^i$ and cluster centroid low rank matrices $\mathbb{L}_j, \mathbb{R}_j$ measured by Eq. (\ref{eq4.4}).
		
	\State Update the centroids $\mathbb{L}_j, \mathbb{R}_j$ for each matrix cluster by Algorithm \ref{alg1} and obtain their corresponding middle low rank matrices $\mathbb{M}^i$ by Eq. (\ref{eq4.7}).
		
	\State Repeat Steps 2-3 until the termination condition is satisfied. 
		
\end{algorithmic}\label{alg3}
\end{algorithm}

Here we check the convergence of Algorithm \ref{alg3} by using the relative reduction of the
the objective value in Eq. (\ref{eq4.8}), i.e. observing if the following inequality holds:

\begin{equation} \label{eq4.9}
	\frac{WCSSRE(t -1) - WCSSRE(t)}{WCSSRE(t -1)} \: \le \: \eta,
\end{equation}

\noindent where $WCSSRE$ stands for the within-cluster sum of square reconstruction error given by

\begin{equation}\label{eq4.10}
	WCSSRE \; := \; \sum_{j=1}^K \sum_{\mathbb{A}^i \in \mathbf{V}_j} \Vert \mathbb{A}^i - \mathbb{L}_j \mathbb{M}^i \mathbb{R}_j^T \Vert_F^2.
\end{equation}

\noindent for some small threshold $\eta > 0$ and $WCSSRE(t)$ denotes the WCSSRE value at the $t$-th iteration of Algorithm \ref{alg3}. We will further study its convergence in Section 4.3.

\subsection{Comparison among CGLRAM, GLRAM and SVD}

We have introduced three data compressing methods in Section 2, the SVD method, the GLRAM method in Algorithm \ref{alg1} and the CGLRAM method in Algorithm \ref{alg3}. Although these three algorithms extract features in different forms, we shall reveal an essential relationship among the SVD, GLRAM and CGLRAM methods. Considering the dataset composed of $N$ matrix patterns $\mathbb{A}^i \in \mathbf{R}^{r \times c}$, we usually convert the matrix $\mathbb{A}^i$ to the vector $a^i = vec(\mathbb{A}^i)$ in the applications to matrix patterns, where $a^i$ denotes the $d$-dimensional ($d=r \cdot c$) vector pattern obtained by sequentially concatenating the columns of $\mathbb{A}^i$. Let $\otimes$ be the Kronecker product, $\Vert \cdot \Vert_F$ and $\Vert \cdot \Vert_2$ denote the Frobenius norm and the Euclidean vector norm respectively.

Generally, SVD computes $k$ orthonormal projection vectors in $\mathbb{P}_{svd} \in \mathbf{R}^{d \times k}$ to extract
$\widetilde{a}_{svd}^i = \mathbb{P}_{svd}^T a^i$ for $a^i$, while $\mathbb{P}_{svd}$ is composed by the first $k$ eigenvectors of $\mathbb{X} = \left[ a^1, ..., a^N \right] $ corresponding to its first largest $k$ eigenvalues. In other words, the reconstructed sample of $a^i$ by SVD can be denoted as $\widetilde{a}_{svd}^i$. Consequently, the projection matrix $\mathbb{P}_{svd}$ is the global optimizer of the following minimization problem, and its objective renders the reconstruction error of the rank-$k$ TSVD:

\begin{equation}\label{eq4.11}
	\begin{aligned}
		&\min \; \; J_{svd}(\mathbb{P}) = \sum_{i=1}^N \Vert a^i - \mathbb{P}\mathbb{P}^Ta^i \Vert_2^2,\\
		&s.t. \quad \; \mathbb{P}^T \mathbb{P} =\mathbb{I}_k,
	\end{aligned}
\end{equation}

\noindent where $\mathbb{I}_k$ denotes the $k \times k$ identity matrix.

Similarly, we could rewrite the optimization structure of GLRAM in Eq. (\ref{eq2.3}) in the vector-pattern form. By Theorem  3.1 in \cite{ye2004generalized}, we extract $\mathbb{M}^i = \mathbb{L}^T \mathbb{A}^i \mathbb{R}$, the GLRAM’s minimization program in Eq. (\ref{eq2.3}) is equivalent to \cite{liu2010generalized}:

\begin{equation}\label{eq4.12}
	\begin{aligned}
		&\min \; \; J_{glram}(\mathbb{P}) = \sum_{i=1}^N \Vert a^i - \mathbb{P}\mathbb{P}^Ta^i \Vert_2^2,\\
		&s.t. \quad \; \mathbb{P}^T \mathbb{P} =\mathbb{I}_k, \; \mathbb{P} = \mathbb{R} \otimes \mathbb{L},\\
		& \qquad \enspace \: \mathbb{L}^T\mathbb{L} = \mathbb{I}_k, \; \mathbb{R}^T \mathbb{R} = \mathbb{I}_k. 
	\end{aligned}
\end{equation}
Then the optimization programming problem of the CGLRAM method in Eq. (\ref{eq4.8}) is written as the following formulation in a similar manner:

\begin{equation}\label{eq4.13}
	\begin{aligned}
		&\min \; \; J_{cglram}(W, \mathbb{P}) = \sum_{j=1}^K \sum_{i=1}^N w_{ij} \Vert a^i - \mathbb{P}_j \mathbb{P}_j^Ta^i \Vert_2^2,\\
		&s.t. \quad \; \mathbb{P}_j^T \mathbb{P}_j =\mathbb{I}_k, \; \mathbb{P}_j = \mathbb{R}_j \otimes \mathbb{L}_j,\\
		& \qquad \enspace \: \mathbb{L}_j^T\mathbb{L}_j = \mathbb{I}_k, \; \mathbb{R}_j^T \mathbb{R}_j = \mathbb{I}_k, \\
		&\qquad \enspace \sum_{j=1}^K w_{ij} = 0, \; i = 1,2,....,N,\\
		&\qquad \enspace \: w_{ij} = 0 \; or \;1,\; i = 1,2,...,N, j = 1,2,...,K,
	\end{aligned}
\end{equation}

\noindent where the decision variable $W = \left[ w_{ij}\right] $ is a $N \times K$ real matrix signing the clustering results 

\begin{equation} \label{eq4.14}
	w_{ij} = \left\{
	\begin{aligned}
		\;1, &\enspace if \;\mathbb{A}^i \; is \; in \; \mathbf{V}_j,\\
		\;0, &\enspace otherwise.
	\end{aligned}
	\right.
\end{equation}

Shown as Eq. (\ref{eq4.13}), we also introduce the decision variable $W$ in Eq. (\ref{eq4.14}) and Eq. (\ref{eq4.2}) to show the clustering result. Then we will obtain different constraint spaces $\mathbf{S}$ for the decision variable $W$. The following conclusions are derived directly from the nature of these three dimensionality reduction techniques. The amount of clusters of GLRAM is regarded as $K=1$, since all matrix patterns share a common projection matrix $\mathbb{P} = \mathbb{R} \otimes \mathbb{L}$, while SVD produces $N$ pairs of low-rank matrices by $\mathbb{A}^i_{svd} = \mathbb{L}^i \mathbb{M}^i \mathbb{R}^i$ for $i = 1,2,...,N$ and we can view as that the dataset is divided into $N$ subsets, where $ \mathbb{L}^i, \mathbb{R}^i$ contain the left and right singular vectors of $\mathbb{A}^i$ respectively. CGLRAM makes $K$ partitions, where $K \le N$ is a specific real number given by us.

\begin{equation} \label{eq4.15}
	\mathbf{S}_{svd} \in \mathbf{R}^{N \times N}, \quad \mathbf{S}_{cglram} \in \mathbf{R}^{N \times K}, \quad \mathbf{S}_{glram} \in \mathbf{R}^{N \times 1}.
\end{equation}

Therefore, we reveal the relationship among CGLRAM, GLRAM and SVD in the following theorem. By Theorem \ref{th4.2}, we can show that GLRAM and SVD achieve the highest and lowest reconstruction errors separately under the same number of reduced dimension, while CGLRAM performs in median. This phenomenon is attributed to the facts that: 1) CGLRAM, GLRAM and SVD in fact optimize the same objective function, and 2) their constraint spaces have similar patterns as shown in Eq. (\ref{eq4.15}) that SVD has the largest constraint space and GLRAM occupies the least, while CGLRAM lies between two parties. 

\begin{theorem}\label{th4.2}
	The objective functions of the SVD, GLRAM and CGLRAM methods are identical while they have been imposed different constraints.
\end{theorem}

To make a better comparison among these three dimensionality reduction techniques, we also summarize the computational and space complexity of the SVD, CGLRAM and GLRAM methods, which is given in Table \ref{tab1}. Although CGLRAM may consume high computational time due to its nested iteration structure, it possesses relatively satisfactory memory storage requirement and reconstruction error than SVD and GLRAM. Therefore, the CGLRAM algorithm achieves a good trade-off among the data compression ratio, the numerical precision and effectiveness in the low rank approximation process.

\begin{table}[!ht] 
	\centering 
	\begin{tabular}{ccc} \hline 
		Method & Time & Space  \\ \hline
		SVD & $O(2MN^3)$ & $Nk^2 +Nk(r+c)$ \\
		CGLRAM & $O(TKMN^2(2I(N+2k)))$ & $Nk^2 + Kk(r+c)$ \\
		GLRAM & $O(MN^2(2I(N+2k)))$ & $Nk^2 + k(r+c)$ \\ \hline
	\end{tabular}
	\caption{Comparison of SVD, CGLRAM and GLRAM: $M, N, k, K, I, T$ denote the amount of samples, dimensions of the previous and low rank data, the number of clusters, and the number of iterations in the while loop in GLRAM and Lloyd's methods respectively.}
\label{tab1}
\end{table}

From Table \ref{tab1}, we could heuristically conclude the reasons why we select to use CGLRAM instead of the standard GLRAM or SVD. First, CGLRAM introduces the concept of clustering into dimensionality reduction, which makes it suitable for more general real-life scenarios. For example, if intermittency is important in data or the dynamics are depicted by less related modes, CGLRAM could still capture essential features and construct a main sketch by imposing a clustering procedure. Second, CGLRAM reduces the amount of computation and saves memory space compared with the full SVD and GLRAM analysis. GLRAM involves the solution of a whole $N\times N$ eigenproblem and SVD makes individual two-sided linear transformations for each sampling matrix, while CGLRAM divides the complex system into several smaller sub-eigenproblems and requires less $\mathbb{L}, \mathbb{R}$s. We also further discuss the numerical performance of these data-compressing algorithms in the next section.

\subsection{Convergence analysis}

In this subsection, we provide some theoretical analysis towards Algorithm \ref{alg3}. We first give an equivalent formulation to the minimization problem in Eq. (\ref{eq4.1}) and then prove the finite convergence of Algorithm \ref{alg3}. We first rewrite the optimization program in the cluster-based formulation. Let $\{\mathbb{A}^i\}_{i=1}^N \in \mathbf{R}^{r \times c}$ be a finite sequence of matrices to be partitioned $K$ clusters with $2<K<N$, then we consider the following mathematical formulation:

\begin{equation}\label{eq4.16}
	\begin{aligned}
		(P): \quad&\min \; f(W,C) = \sum_{j=1}^K \sum_{i=1}^{N} w_{ij}D_{ij},\\
		&s.t. \: \sum_{j=1}^K w_{ij} = 0, \; i = 1,2,....,N,\\
		&\qquad w_{ij} = 0 \; or \;1,\; i = 1,2,...,N, j = 1,2,...,K,
	\end{aligned}
\end{equation}

\noindent where the decision variable $W = \left[ w_{ij}\right] $ is a $N \times K$ real matrix defined in Eq. (\ref{eq4.11}) and $C =  \left[ C_1,C_2,...,C_K\right] $ denotes the centroid of clusters with $C_j = \{\mathbb{L}_j, \mathbb{R}_j\}$. The distance between sample matrix $\mathbb{A}^i$ and cluster $\mathbf{V}_j$ is defined as Eq. (\ref{eq4.5}).

To investigate the properties of Program (P), we start with the following definition towards its objective function. Definition \ref{def4.1} reveals that Program (P) can be reduced as a univariate programming problem with respect to $W$. Actually, by Theorem \ref{th4.1}, if we fix $W$, the objective function can be transformed into the sum of several GLRAM goals, i.e., the decision variable $C$ can be determined by $W$. We also detect the relationship between $W$ and $C$ in Theorem \ref{th4.4}.

\begin{definition} \label{def4.1}
	The reduced function $F(W)$ of Program (P) is defined as follows and $W$ is a $N \times K$ real matrix
	
	$$F(W) \;= \; \min_{C \in \Omega} \left\lbrace f(W, C) \right\rbrace ,$$
	
	\noindent where $\Omega$ is the constraint space defined by $\mathbb{L}_j \in \mathbf{R}^{r \times k}, \mathbb{L}_j^T \mathbb{L}_j = I_{k}$ and $\mathbb{R}_j \in \mathbf{R}^{k \times c}, \mathbb{R}_j^T \mathbb{R}_j = I_{k}$ for all $j = 1,2,...,K$.
	
\end{definition}

Based on \cite{selim1984k}, we can also conclude the property of the constraint space of Program (P).

\begin{theorem} \label{th4.3}
	Consider the set $\mathbf{S}$ given by
	
	\begin{equation}
		\begin{aligned}
			\mathbf{S} = \left\{  \sum_{j=1}^K w_{ij} = 1, i = 1,2,...,N ; \quad w_{ij} \ge 0, i = 1,2,...,N, j = 1,2,...,K \right\} .
		\end{aligned}
		\nonumber
	\end{equation}
	
	\noindent The extreme points of $\mathbf{S}$ satisfy constraints in Program (P).
\end{theorem}

Since the function $F$ contains all properties of $f$ and also there exists an extreme point solution of Program (RP) which in turn satisfies constraints in Program (P). We can immediately conclude the following important optimization problem that is equivalent to Program (P), and any results in the rest of the section corresponding to one of the problems can be transformed into the other.

\begin{theorem} \label{th4.4}
	
	The reduced problem (RP) defined as follows and Program (P) are equivalent.	
	
	\begin{equation}\label{eq4.17}
		(RP): \quad\min \; F(W), \enspace subject \; to \;W \in \mathbf{S}. 
	\end{equation}
\end{theorem}

In the next process, we characterize partial optimal solutions of Program (P) and show their relationship with the Kuhn-Tucker points.

\begin{definition} [\cite{lasdon2002optimization}] \label{def4.2} 
	A point $(W^*, C^*)$ is a partial optimal solution of Program (P) if it satisfies that
	
	\begin{equation} \label{eq4.15}
		\begin{aligned}
			&f(W^*,C^*) \le f(W,C^*), for \; all \;W\in \mathbf{S},\\
			&f(W^*,C^*) \le f(W^*,C), for\; all \;C\in \Omega,
		\end{aligned}
	\end{equation}
	
\end{definition}

By the idea from \cite{selim1984k}, the following theorem could reveal the connection between these partial optimal solutions $W^*, C^*$ and the Kuhn-Tucker points of Program (P).

\begin{theorem}\label{th4.5}
	Suppose $F(W^*, C)$ is differentiable at $C = C^*$, then $(W^*,C^*) $ is a Kuhn-Tucker point of Program (P) if and only if it is a partial optimal solution of Program (P).
\end{theorem}

\begin{proof}
	First assume $(W^*, C^*)$ is a partial optimal solution satisfying the partially differentiability condition. Then $(W^*, C^*)$ is obtained by the following two problems:
	
	\begin{equation}\label{eq4.18}
		\begin{aligned}
			& (P_1): \quad Given \; \hat{C} \in \Omega, \; \min \; f(W,\hat{C}), \enspace subject \; to \;W \in \mathbf{S}. \\
			& (P_2): \quad Given \; \hat{W} \in \mathbf{S}, \; \min \; f(\hat{W},C), \enspace subject \; to \;C \in \Omega. 
		\end{aligned}
	\end{equation}
	
	Since $W^*$ solves ($P_1$) for $\hat{C} = C^*$, then $W^*$ satisfies Kuhn-Tucker conditions of ($P_1$). Similarly, if $C^*$ solves ($P_2$), $C^*$ satisfies its Kuhn-Tucker conditions when $\hat{W} = W^*$. Then it is obvious that $(W^*, C^*)$ is a Kuhn-Tucker point of Program (P).
	
	Then let  $(W^*, C^*)$ be a Kuhn-Tucker point of Program (P), and suppose it is not a partial optimal solution of (P). Then $W^*, C^*$ cannot solve ($P_1$) and ($P_2$) respectively, i.e. they do not satisfy the corresponding Kuhn-Tucker condition, which causes a contradiction. Therefore, $(W^*, C^*)$ is indeed a partial optimal solution of (P).
	
\end{proof}

Finally, we come to show our cluster-based data compressing algorithm has the finite convergence.

\begin{theorem}\label{th4.6}
	The CGLRAM method summarized in Algorithm \ref{alg3} converges to a partial optimal solution of Program (P) in a finite number of iterations.
\end{theorem}

\begin{proof}
	First suppose the iterative algorithm proceeds from iteration $t$ to iteration $t+1$, then it suffices to hold that $f(W^{t+1}, C^{t+1}) \le f(W^{t}, C^{t})$. On the one hand, the inequality of  $f(W^{t+1}, C^{t}) \le f(W^{t}, C^{t})$ follows immediately by the logic of Algorithm \ref{alg3}, since $C$ should be updated if there exists sample matrix which finds closer cluster than the cluster previously assigned to it and we make different tessellation from $W^{t}$ to $W^{t+1}$. On the other hand,  $f(W^{t}, C^{t+1}) \le f(W^{t}, C^{t})$ is directly obtained by Theorem \ref{th2.3} that Algorithm \ref{alg1} monotonically non-increases the RMSRE value, hence it converges in the limit.
	
	Meanwhile, each extreme point of $\mathbf{S}$ will be visited at most once before the termination. Therefore, Algorithm \ref{alg3} will reach a partial optimal solution in a finite number of iterations, since there are a limited quantity of extreme points of $\mathbf{S}$. Based on the definitions of partial optimal solutions depicted in Eq. (\ref{eq4.15}), the final clustering result is locally optimal.
	
\end{proof}

\section{Numerical experiments}

In this section, we experimentally evaluate our CGLRAM algorithm in Algorithm \ref{alg3}. All of our simulations are carried out by using MATLAB R2022a software on an Apple M1 machine with 8GB of memory. To thoroughly illustrate CGLRAM’s performance, we consider the following two specific applications: image compressing of a hand-written image dataset and dimensionality reduction on snapshot matrices of the stochastic Navier-Stokes equation.

\subsection{Applications on image compressing}

We first present a detailed comparative study on a real-world image dataset among three data compressing algorithms including the standard GLRAM, the proposed CGLRAM and K-means + GLRAM algorithm, where the pre-classification is conducted by K-means clustering. The comparison focuses on the reconstruction error (measured by WCSSRE) and clustering quality (measured by error reduction rate).

\begin{figure}[ht]
	\centering
	\includegraphics[scale=0.5]{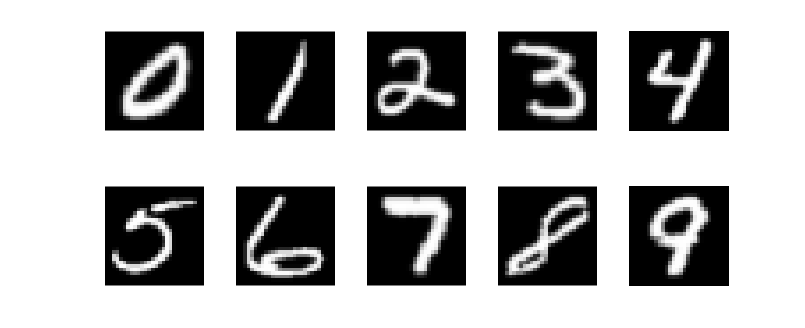}
	\caption{Handwritten figure examples in EMNIST-digits.}
	\label{fig5.1}
\end{figure}

The image dataset EMNIST \cite{cohen2017emnist} is adopted in the experiment. It is a set of handwritten character digits derived from the NIST Special Database 19 and converted to a 28x28 pixel image format, which contains digits, uppercase and lowercase handwritten letters. The EMNIST dataset takes advantage of its understandable and intuitive structure, relatively small size and storage requirements, accessibility and ease-of-use property, which makes it become a standard benchmark for learning, classification and computer vision systems in recent years. In our simulation, we randomly select 1000 samples from its sub-dataset \emph{emnist-digits.mat} and the test collection includes handwritten images of the numbers 0-9. We summarize the statistics of our test dataset and visualize its basic classes in Table \ref{tab5.1} and Figure \ref{fig5.1} respectively.

\begin{table}[ht]
\centering 
\begin{tabular}{cccc}
\hline
	Dataset & Size ($N$) & 	Dimension ($r \times c$) & Number of classes ($K$)\\ \hline
    EMNIST-digits & 1000 & 28 $\times$ 28 & 10\\ \hline
\end{tabular}
\caption{Statistics of our test dataset 1.}
\label{tab5.1}
\end{table}

\begin{figure}[ht]
	\centering
	\includegraphics[scale=0.33]{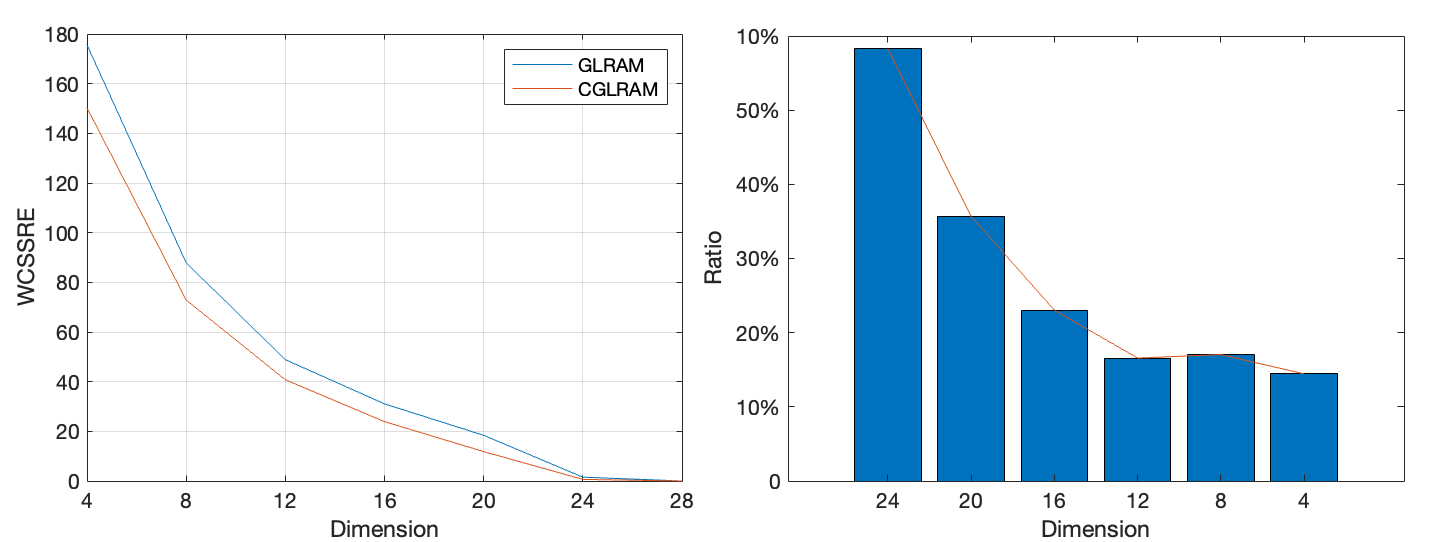}
	\caption{Comparison of reconstruction error (left) and error reduction rate (right) between GLRAM and CGLRAM.}
	\label{fig5.2}
\end{figure}

We first evaluate the effectiveness of the proposed CGLRAM algorithm in terms of the reconstruction error measured by WCSSRE and compare it with the standard GLRAM. Here both GLRAM and CGLRAM employ the same numbers of reduced dimensions$\{4,8,12,16,20,24,28\}$. The comparative results are concluded in Figure \ref{fig5.2}, where the horizontal axis denotes the amount of reduced dimension (k), and the vertical axis denotes the WCSSRE value (left) and the rate of decline in WCSSRE of these two methods (right). It is observed that CGLRAM has a lower WCSSRE value in all cases than GLRAM, which corresponds to our conclusion in Section 4.2. We can also achieve a very high enhancement rate in matrix reconstruction accuracy, nearly up to 60\% when adopting $k = 24$. In other words, our proposed method is more competitive with the standard GLRAM in image compression. 

\begin{figure}[ht]
\centering
\includegraphics[scale=0.33]{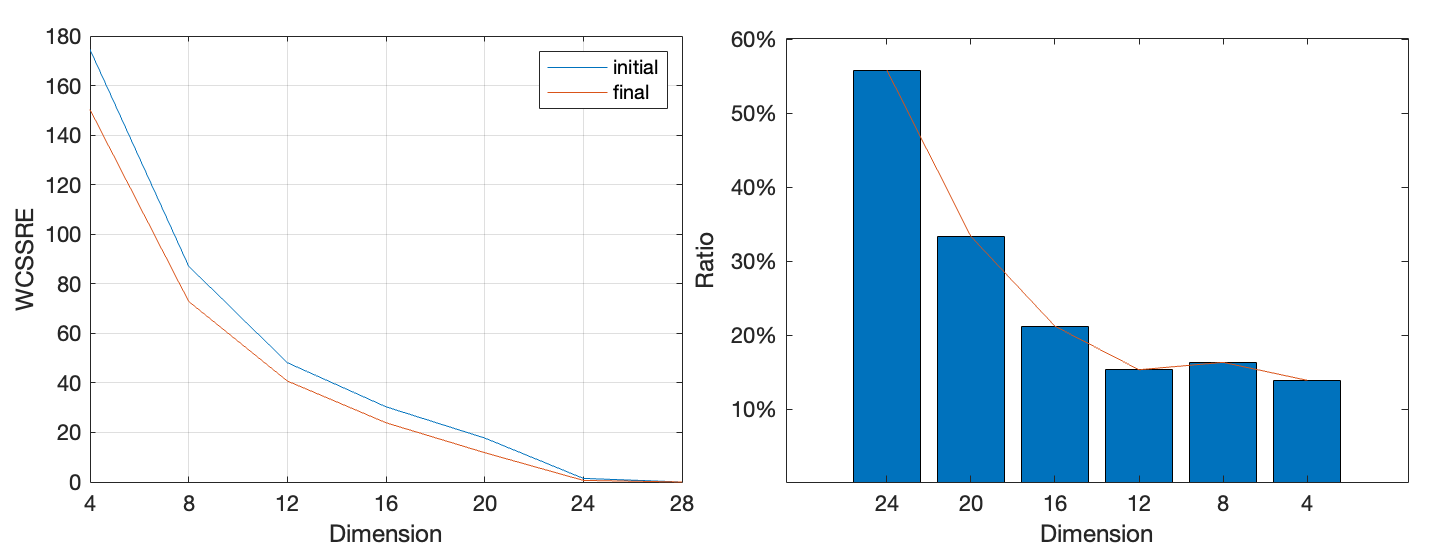}
\caption{Comparison of initial and final CGLRAM objective values (left) and its enhancement rate (right).}
\label{fig5.3}
\end{figure}

In our simulations, the effectiveness of pre-clustering procedure in CGLRAM is also studied from two aspects. First we compare the initial and final CGLRAM objective values (WCSSRE) and compute their corresponding ratios as illustrated in Figure \ref{fig5.3}. We can observe that CGLRAM truly updates better classifications and significantly reduces matrix reconstruction errors with the help of clustering process. Second, a comparative study is also provided between CGLRAM and K-means+GLRAM, where the images are pre-classified by the traditional K-means clustering method and then compressed by GLRAM in each cluster. Figure \ref{fig5.3} presents that CGLRAM performs better in image clustering and it levels up matrix reconstruction accuracy compared with K-means+GLRAM. This may be related to the fact that CGLRAM is able to utilize correlation information intrinsic in the image collection, and it generalizes standard K-means to a broader space and similarity metric accordingly, which leads to better classification performance.

\begin{figure}[ht]
	\centering
	\includegraphics[scale=0.33]{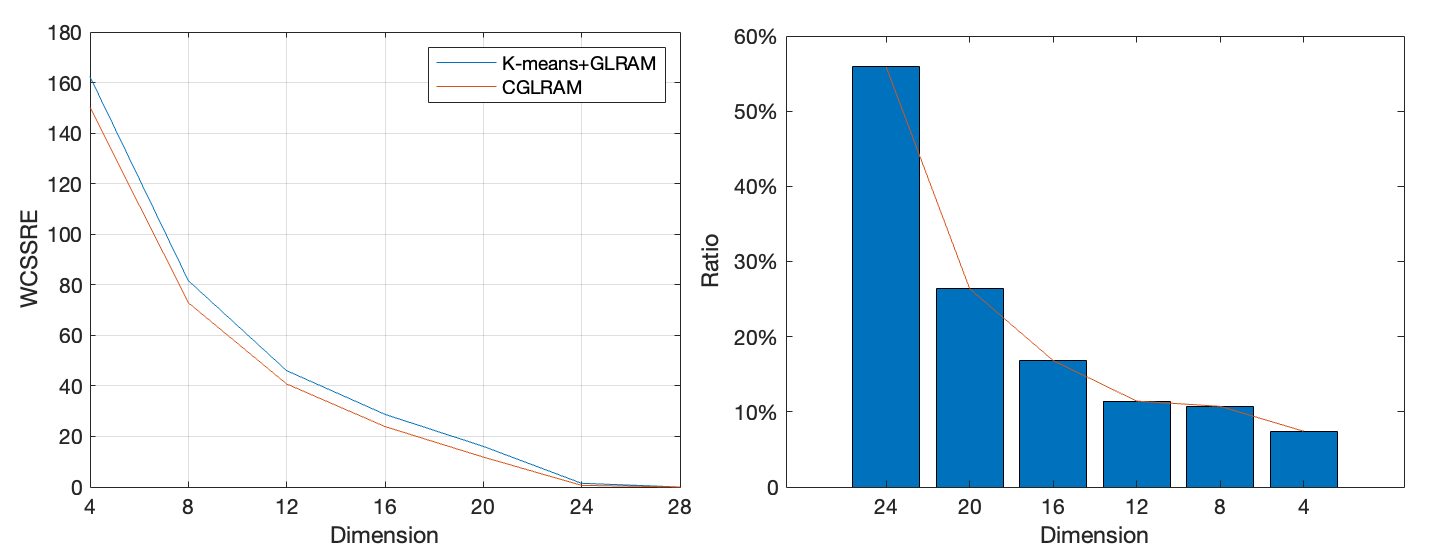}
	\caption{Comparison of reconstruction error between K-means+GLRAM and CGLRAM.}
	\label{fig5.4}
\end{figure}

\begin{figure}[ht]
	\centering
	\includegraphics[scale=0.45]{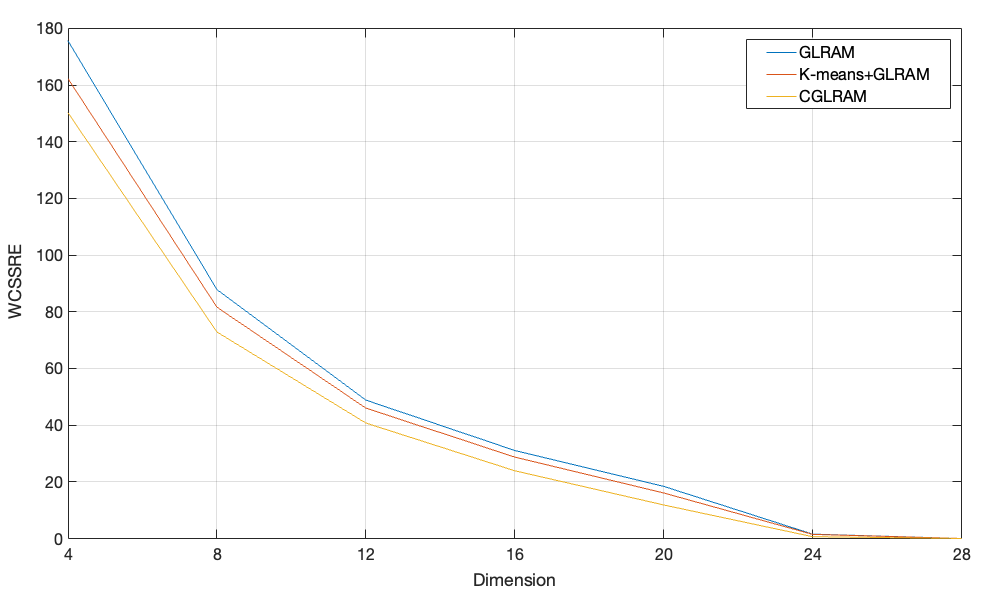}
	\caption{Comparison of reconstruction error among GLRAM, K-means+GLRAM and CGLRAM.}
	\label{fig5.5}
\end{figure}

Finally, we summarize the simulation results of these three methods in Figure \ref{fig5.5} and Table \ref{tab5.2}, which demonstrates the within-cluster sum of square reconstruction errors and corresponding error reduction rates obtained by $\displaystyle {(err_{initial} - err_{final})}/{err_{initial}}$ of three different compression methods, GLRAM, K-means+GLRAM and CGLRAM. It is observed that CGLRAM can achieve the highest dimension reduction accuracy and significantly decreases matrix reconstruction errors than other traditional methods, which validates the effectiveness of our proposed algorithm. The numerical results also show that matrix reconstruction errors decrease monotonically for all methods with the growth of reduced dimension $k$. Generally, a large $k$ could improve the performance of CGLRAM in reconstruction and classification, while its computation cost also increases as $d$ increases. Therefore we need to consider the tradeoff between the computational precision and complexity in choosing the best $d$.

\begin{table}[ht]
	\centering 
	\resizebox{\textwidth}{!}{
	\begin{tabular}{c|c|ccccccc}
		\hline
		\multicolumn{2}{c|}{Dimension $k$} & $28$ & $24$ & $20$ & $16$ & $12$ & $8$ & $4$ \\ \hline
		 \multirow{3}{*}{WCSSRE} & GLRAM & $4.5045 \times 10^{-13}$ & $1.5443$ &	$18.4009$ &	$31.0613$ & $48.8563$ & $87.8599$ &	$175.6365$\\ 
		 & K-means+GLRAM & $4.3623 \times 10^{-13}$ & $1.4606$ & $16.0694$ &	$28.7381$ & $46.0196$ & $81.6403$ &	$162.2953$\\ 
		 & CGLRAM & $4.2015 \times 10^{-13}$ & $0.6442$ & $11.8328$ &	$23.9073$ & $40.7598$ & $72.8806$ &	$150.2670$\\  \hline
		  \multirow{3}{*}{Error reduction ratio} &\multicolumn{2}{c|}{GLRAM vs CGLRAM} & $58.2853\%$ & $35.6941\%$ & $23.0320\%$ &	$16.5721\%$ & $17.0491\%$ & $14.4443\%$ \\ 
		  &\multicolumn{2}{c|}{K-means+GLRAM vs CGLRAM} & $55.8952\%$ & $26.3639\%$ & $16.8099\%$ &	$11.4296\%$ & $10.7296\%$ & $7.4114\%$ \\ 
		  &\multicolumn{2}{c|}{initial CGLRAM vs final CGLRAM} & $55.8107\%$ & $33.3608\%$ & $21.2228\%$ &	$15.3444\%$ & $16.3450\%$ & $13.9174\%$ \\ 	\hline
	\end{tabular}}
	\caption{Simulation results for 3 different dimension reduction methods GLRAM, K-means+GLRAM and CGLRAM about within-cluster sum of square reconstruction errors (WCSSRE) and their error reduction rates obtained by $(err_{ini} - err_{fin})/err_{ini}$.}
	\label{tab5.2}
\end{table}

\subsection{Applications on stochastic partial differential equations (SPDE)}

Another numerical simulation is performed for chaotic complex systems, where their solutions concern with the physical, temporal and random spaces and tend to be represented as large and complex matrices in real implementation. Let $D$ be a bounded open set, we consider the following two-dimensional stochastic Navier-Stokes equation with additive random forcing:

\begin{equation}
\begin{aligned}\label{eq5.1}
&\bm{u}_t - \nu \Delta \bm{u} + (\bm{u} \cdot \nabla)\bm{u} + \nabla p = \bm{f} + \sigma \dot{\bm{W}},\\
&\nabla \cdot \bm{u} = 0, 
\end{aligned}
\end{equation}

\noindent with initial velocity

\begin{equation}
\bm{u}\vert_{t=0} = \bm{u}_0,
\nonumber
\end{equation}

\noindent and Dirichlet boundary conditions

\begin{equation}
\bm{u}\vert_{\partial D \times (0,T]} = \bm{0},
\nonumber
\end{equation}

\noindent where $\nu > 0$ denotes the viscosity parameter, $\bm{u} = [u, v]^T$, $\bm{u}_0 = [\phi_1, \phi_2]^T$, $\bm{x} = [x, y]^T$, $\sigma = diag[\sigma_1(\bm{x}), \sigma_2(\bm{x})]$, and $\bm{W}$ is a Brownian motion vector. Hence the term $\dot{\bm{W}}$ represents the additive random term, and $\bm{f} = [f_1, f_2]^T$ is the external force. For the representation of the white noise $\dot{\bm{W}}(\bm{x}, t)$, we employ the piecewise constant approximation \cite{du2002numerical} defined as

\begin{equation} \label{eq5.2}
	\dot{\bm{W}}(\bm{x}, t) = \frac{1}{\sqrt{\Delta t}} \sum_{n = 0}^{N-1} \chi_n(\bm{x}, t) \eta_n(\omega),
\end{equation}

\noindent where, for $n = 0, 1,\cdots,N - 1$, the independent and identically distributed (i.i.d.)
random variable $\eta_n$ satisfies the standard normal distribution $N(0, 1)$ and the characteristic function $\chi_n$ is given by 

\begin{equation}
	\chi_n(t) = \left\{
	\begin{aligned}
		& 1 \quad \text{if} \enspace t \in [t_n, t_{n+1}),\\
		& 0 \quad \text{otherwise}.
	\end{aligned}
	\right.
	\nonumber
\end{equation}

We first provide some required definitions of function spaces and notations, and define the stochastic Sobolev space as

\begin{equation} \label{eq5.3}
	\mathcal{L}^2(\Omega; 0, T; H^1(D)) = L^2(0, T; H^1(D)) \otimes L^2(\Omega),
\end{equation}

\noindent which is a tensor product space, here $\Omega$ is the stochastic variable space, and $L^2(0, T; H^1(D))$ is the space of strongly measurable maps $\bm{v} : [0, T] → H^1(D)$ such that

\begin{equation} \label{eq5.4}
	\Vert v \Vert^2_{L^2(0, T; H^1(D))} = \int_0^T \Vert v \Vert^2_D \rm{d}t,
\end{equation}

\noindent so the norm of the tensor Sobolev space in Eq. (\ref{eq5.3}) can be defined as

\begin{equation} \label{eq5.5}
	\Vert v \Vert^2_{\mathcal{L}^2(\Omega; 0, T; H^1(D))} = \int_0^T \Vert v \Vert^2_{L^2(0, T; H^1(D))} \rm{d}P,
\end{equation}

Then, we can obtain the following abstract problem of the stochastic NS equations in Eq. (\ref{eq5.1})

\begin{equation} \label{eq5.6}
	\rm{d} \bm{u} = \mathcal{A}(\bm{u}) \rm{d}t + \mathcal{B}(\bm{x}) \rm{d}W,
\end{equation}

\noindent by applying operators 

\begin{equation}
\begin{aligned}
& \mathcal{A}(\bm{u}) =  \nu \Delta \bm{u} - (\bm{u} \cdot \nabla)\bm{u} - \nabla p + \bm{f},\\
& \mathcal{B}(\bm{x}) = \sigma(\bm{x}).
\nonumber
\end{aligned}
\end{equation}

To discretize Eq. (\ref{eq5.6}) with regards to time, we divide the time interval $[0,T]$ into $L$ subintervals of duration $\Delta t = T / L$. Thus for $l = 0,\cdots,L-1$, given $\bm{u}^{l} = \bm{u}(t_l, \bm{x})$ and $t_{l+1} = t_l + \Delta t$, the solution of the stochastic differential equation system is provided by

\begin{equation} \label{eq5.7}
 \bm{u}^{l+1} =  \bm{u}^{l} + \int_{t_l}^{t_{l+1}} \mathcal{A}(\bm{u}(s)) \rm{d}s + \mathcal{B}(\bm{x}) \int_{t_l}^{t_{l+1}} \rm{d}W.
\end{equation}

\noindent Here we choose the order 1 weak implicit Euler scheme \cite{kloeden1992stochastic} to tackle with the integral term above

\begin{equation} \label{eq5.8}
	\bm{u}^{l+1} =  \bm{u}^{l} + \frac{\Delta t}{2} \left( \mathcal{A}(\bm{u}^{l+1}) + \mathcal{A}(\bm{u}^{l})\right) + \mathcal{B}(\bm{x}) \Delta W.
\end{equation}

Simultaneously, the standard finite element method is adopted for spatial discretization and we linearize the nonlinear convective term by Newton’s method \cite{gunzburger1996navier}. Then, a weak formulation of the stochastic NS equations in Eq. (\ref{eq5.1}) is given as follows: for $l = 0,\cdots,L-1$, determine the pairs $(\bm{u}^{l+1}, p^{l+1})$ by solving the linear system 

\begin{equation} \label{eq5.9}
\begin{aligned}
	& (\bm{u}^{l+1} - \bm{u}^{l+1}, \bm{v}) + \frac{\Delta t}{2} \left( a(\bm{u}^{l+1}, \bm{v}) + a(\bm{u}^{l}, \bm{v}) \right) + \frac{\Delta t}{2} \left( c(\bm{u}^{l+1}, \bm{u}^{l}, \bm{v}) + c(\bm{u}^{l}, \bm{u}^{l+1}, \bm{v}) \right) \\
	& \qquad \quad + \frac{\Delta t}{2} \left( b(\bm{v}, p^{l+1}) + b(\bm{v}, p^{l}) \right) = \frac{\Delta t}{2} \left( (\bm{f}^{l+1}, \bm{v}) + (\bm{f}^{l}, \bm{v}) \right) + (\sigma \Delta W, \bm{v}),\\
	& b(\bm{u}^{l+1}, q) = 0.
\end{aligned}
\end{equation}

\noindent for all test functions $\bm{v} \in \mathcal{L}_0^2(\Omega; 0, T; H^1(D))$ and $q \in \mathcal{L}^2(\Omega; 0, T; L^2(D)$, denote $\bm{f}^{l} = \bm{f}(t_l)$, and where we define the bilinear operators

\begin{equation}
	a(\bm{u}, \bm{v}) = \int_D \nu \nabla \bm{u} \nabla \bm{v} dD,
	\nonumber
\end{equation}

\noindent and

\begin{equation}
	b(\bm{v}, q) = - \int_D q \nabla \bm{v} dD,
	\nonumber
\end{equation}

\noindent and the trilinear operator

\begin{equation}
	c(\bm{w}, \bm{u}, \bm{v}) = \int_D \bm{w} \nabla \bm{u} \bm{v} dD,
	\nonumber
\end{equation}

\noindent and the $\mathcal{L}^2(D)$ inner product

\begin{equation}
	(\bm{u}, \bm{v}) = \int_D \bm{u} \bm{v} dD.
	\nonumber
\end{equation}

To check the validation of our methodology, a numerical example is tested with the following settings. Let the domain be $D = [0,1]^2$ and we set $T = 3$, $\Delta t = 0.01$, $\nu = 0.01$,  and the number of nodes in our finite element mesh is $\# \text{ of FEM nodes} = 665$. The initial conditions are given by

\begin{equation}
\phi_1(x, y) = - \varphi(x) \varphi'(y), \quad \phi_2(x, y) = \varphi'(x) \varphi(y),
\nonumber
\end{equation}

\noindent where $\varphi(z) = 10 z^2 (1 - z)^2$, and let $\sigma$ be defined as

\begin{equation}
\sigma_1(x, y) = \cos(x) \sin(y), \quad \phi_2(x, y) = \sin(x) \sin(y).
\nonumber
\end{equation}


The number of samples in our simulation is $N = 100$, and the statistics of this test dataset are provided in Table \ref{tab5.3}. Here we do not apply a very high sample size $N$, since sample matrices have a relatively large dimension. Each sample represents a Navier-Stokes numerical solution $\bm{u} = (u, v)^T$, while $u,v \in \mathbf{R}^{665 \times 331}$, and therefore the overall data amount is still large.

\begin{table}[ht]
\centering 
\begin{tabular}{cccc}
\hline
Dataset & Size ($N$) & 	Dimension ($r \times c$) & Number of classes ($K$)\\ \hline
NS snapshots $u, v$ & 100 & 665 $\times$ 301 & 2\\ \hline
\end{tabular}
\caption{Statistics of our test dataset 2.}
\label{tab5.3}
\end{table}

\begin{figure}[ht]
\centering
\includegraphics[width=1\textwidth]{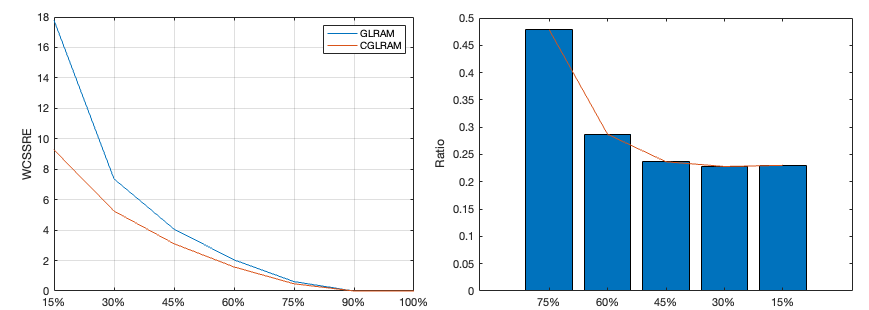}
\caption{Comparison of reconstruction error (left) and error reduction rate (right) between GLRAM and CGLRAM under 7 different dimension reduction ratios $\tau = \{100\%, 90\%, 75\%, 60\%, 45\%, 30\%, 15\%\}$.}
\label{fig5.7}
\end{figure}

We first compare the simulation results of standard GLRAM and our proposed algorithm in this SPDE setting and the results are demonstrated as follows. The left plot in Figure \ref{fig5.7} visualizes the matrix reconstruction errors obtained by GLRAM and CGLRAM respectively with 7 different dimension reduction ratios $\tau =$ 100\%, 90\%, 75\%, 60\%, 45\%, 30\%, 15\%, and the right plot depicts the corresponding error reduction rates. We observe that the blue line is above the red line in any value of $\tau$ and error reduction ratios between these two methods are also notable around $(20\%, 50\%)$, which presents that CGLRAM has significantly better numerical performances than GLRAM. Table \ref{tab5.3} shows the specific values of WCSSRE and error reduction rates obtained by GLRAM and CGLRAM under different dimension reduction ratios $\tau$ and reduced dimension $k$. Here we do not provide the dimension reduction rates with $\tau = 100\%$ and $\tau = 75\%$, because both GLRAM and CGLRAM errors are so tiny and close that the rates cannot convey useful information in these cases. To sum up, compared with GLRAM, our proposed algorithm CGLRAM significantly reduces reconstruction errors and makes a better balance between dimension reduction rate and computational accuracy. Moreover, the simulation results illustrate the sensitivity of CGLRAM to the choice of $\tau$ and the trend is observed that WCSSRE drops monotonically with the rise of $\tau$, which also corresponds to our findings in Section 5.1.

\begin{table}[ht]
\centering 
\resizebox{\textwidth}{!}{
\begin{tabular}{c|ccccccc}
\hline
dimension reduction ratio $\tau$ & $100\%$ & $90\%$ & $75\%$ & $60\%$ & $45\%$ & $30\%$ & $15\%$ \\ \hline
reduced dimension $k$ & $665$ & $599$ &	$499$ &	$399$ & $299$ & $200$ &	$100$\\ 
WCSSRE of GLRAM & $1.6686 \times 10^{-11}$ & $1.3101 \times 10^{-6}$ & $0.6167$ &	$2.0476$ & $4.0693$ & $7.3672$ & $17.7917$\\ 
WCSSRE of CGLRAM & $1.2948 \times 10^{-13}$ & $1.7022 \times 10^{-10}$ & $0.4748$ &	$1.5806$ & $3.1060$ & $5.2527$ & $9.2751$\\  \hline
\multicolumn{3}{c|}{Error reduction ratio of GLRAM vs CGLRAM} & $47.8685\%$ & $28.7017\%$ & $23.6705\%$ &	$22.8047\%$ & $23.0114\%$ \\ 	\hline
\end{tabular}}
\caption{Simulation results for GLRAM and CGLRAM about within-cluster sum of square reconstruction errors (WCSSRE) and their error reduction rates obtained by $(err_{ini} - err_{fin})/err_{ini}$ under different dimension reduction ratios $\tau$ and reduced dimension $k$.}	\label{tab5.3}
\end{table}

\begin{figure}[!h]
\centering
\includegraphics[width=1\textwidth]{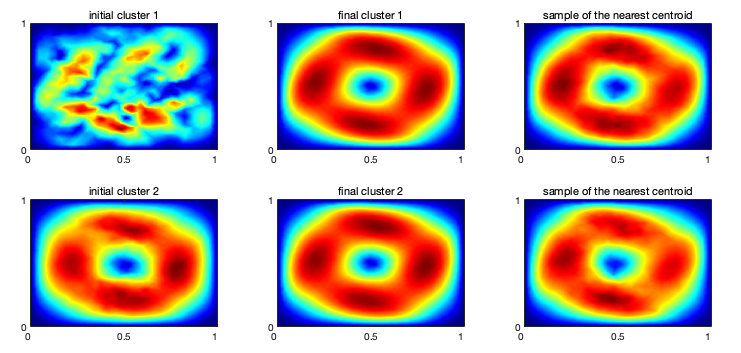}
\caption{The initial cluster mean (left), the final cluster mean (middle), and a random sample in each cluster (right).}
\label{fig5.8}
\end{figure}

\begin{figure}[!ht]
\centering
\includegraphics[width=1\textwidth]{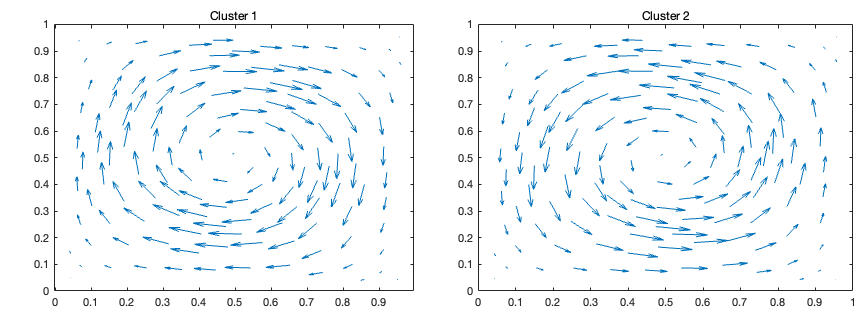}
\caption{The mean of stochastic Navier-Stokes flows of Cluster 1 (left) and Cluster 2 (right).}
\label{fig5.9}
\end{figure}

Similarly, we also study the clustering performance of CGLRAM in this numerical test and demonstrate the results in graphical form. Figure \ref{fig5.8} compares the initial and final classifications of the test dataset and presents a random sample in each cluster. It is obvious that the initial and final centroids are quite different and our clustering procedure is effective in classifying the dataset. Figure \ref{fig5.9} presents the numerical simulation of the mean of the velocity fields in the square domain from two different clusters. It is clear that two subplots in Figure \ref{fig5.9} show a vortex-shape flow, while the flow directions are opposite. The left plot has a tendency to be clockwise and the right plot turns counterclockwise. In other words, the final two clusters own unique and different features and therefore it validates the effectiveness of our clustering procedure. 

In this numerical test, we do not provide the comparative statistics between K-means+GLRAM and CGLRAM, because we find that they have similar classifying and compressing effects. In the previous image experiment in Section 5.1, CGLRAM significantly outperforms K-means+GLRAM and standard GLRAM, while the advantage is not obvious in this test. We think it may result from the fact that SPDE snapshots tend to be large sparse matrices, whose entries have small values. Therefore, discrepancies among these sample matrices become less obvious especially when we use low reduced dimensions. Meanwhile, compared with SPDE snapshots, CGLRAM can better utilize locality information intrinsic in images, for example, smoothness in images, which also leads to superior CGLRAM performance in the previous test.

Another noteworthy issue is that we need to focus on the initial classification in matrix collections since our proposed  CGLRAM  method is sensitive to the initial classification. Unsuitable choice of initial clusters may result in a time-consuming iteration or even stuck in a local optimizer, which is one of the drawbacks of Lloyd's method. Therefore, how to select suitable initial clusters is an essential topic to be explored in our future work.

\section{Conclusions}

In this study, we present a clustering-based generalized low-rank approximation method that introduces the idea of clustering and extends the conventional GLRAM framework. The novel low rank approximation method employs a generalized form of centroids and similarity measures in clustering, thereby facilitating its application to more types of datasets. Relatively to GLRAM, CGLRAM inherits better computational efficiency and precision by pre-classifying the datasets and enhancing correlations within clusters. In summary, CGLRAM leverages advantages from both GLRAM and K-means methods, and presents nice numerical performance in our experiments.

Realistic simulations of complex systems governed by partial differential equations must account for uncertain features of the modeling phenomena. Naturally, reduced-order models are considered to lessen the computational cost of stochastic partial differential equations (SPDEs) by using few degrees of freedom. The objectives here are about the introduction of our low rank approximation method in solving nonlinear or unsteady SPDEs, so that the large and complex matrix calculations can be greatly alleviated. The issues discussed in Section 5 will be further studied in the future.

\section*{Acknowledgments}
The authors would like to thank the anonymous referees and the editor for their valuable comments and suggestions, which led to considerable improvement of the article.

\section*{Conflict of Interest}
All authors declare that there are no conflicts of interest regarding the publication of this paper.

\bibliographystyle{elsarticle-harv}
\bibliography{bib}

\end{document}